# Moments of convex distribution functions and completely alternating sequences

## Alexander Gnedin[1] and Jim Pitman[*,2]

*Mathematisch Instituut and University of California, Berkeley*

**Abstract:** We solve the moment problem for convex distribution functions on $[0,1]$ in terms of completely alternating sequences. This complements a recent solution of this problem by Diaconis and Freedman, and relates this work to the Lévy-Khintchine formula for the Laplace transform of a subordinator, and to regenerative composition structures.

## 1. Introduction

It is well known that the distribution function $F$ of a probability measure on the unit interval $[0,1]$ is uniquely determined by its sequence of moments

$$(1.1) \qquad c(n) := \int_{[0,1]} x^n \, dF(x) , \quad n = 0, 1, \ldots$$

defined by Lebesgue-Stieltjes integration. A complete characterization of such moment sequences was discovered by Hausdorff [15]. To recall his result, for $j = 0, 1, 2, \ldots$, let $\nabla^j$ denote the $j$th iterate of the difference operator

$$\nabla c(n) := c(n) - c(n+1), \quad n = 0, 1, 2, \ldots.$$

That is,

$$(1.2) \qquad \nabla^j c(n) = \sum_{i=0}^{j} (-1)^i \binom{j}{i} c(n+i), \quad n \geq 0, j \geq 0.$$

**Definition 1.1.** A sequence $c = (c(0), c(1), \ldots)$ is *completely monotone* if

$$(1.3) \qquad \nabla^j c(n) \geq 0 \quad \text{for all } j \geq 0, n \geq 0.$$

**Theorem 1.2** (Hausdorff [15]). *A sequence $c$ can be represented as the moment sequence* (1.1) *for some probability distribution $F$ on $[0,1]$ if and only if $c(0) = 1$ and $c$ is completely monotone.*

---

*Supported in part by NSF Grant DMS-04-05779.
[1]Mathematisch Instituut, Rijksuniversiteit Utrecht, P.O. Box 80 010, 3508 TA Utrecht, The Netherlands, e-mail: gnedin@math.uu.nl
[2]University of California, Department of Statistics, 367 Evans Hall #3860, Berkeley, CA 94720-3860, USA, e-mail: pitman@stat.berkeley.edu

A sequence $c$ generates a *triangular array* $(c(n,m);\ 0 \le m \le n,\ n = 0, 1, \ldots)$ by the formula

$$(1.4) \qquad c(n,m) := \binom{n}{m} \nabla^{n-m} c(m), \qquad 0 \le m \le n.$$

Obviously, $c$ is completely monotone if and only if

$$(1.5) \qquad c(n,m) \ge 0, \qquad 0 \le m \le n.$$

For such a sequence $c$ with $c(0) = 1$, for each $n = 0, 1, \ldots$ the numbers $(c(n,m), 0 \le m \le n)$ form the probability distribution of a random variable $S_n$ with values $0, 1, \ldots, n$:

$$(1.6) \qquad c(n,m) = \mathbb{P}(S_n = m), \qquad 0 \le m \le n.$$

Hausdorff showed that then there is the convergence in distribution,

$$(1.7) \qquad \lim_{n \to \infty} \mathbb{P}(S_n/n \le x) = F(x),$$

at all continuity points $x$ of the unique probability distribution function $F$ whose $n$th moment is $c(n)$. Moreover,

$$(1.8) \qquad c(n,m) = \binom{n}{m} \int_{[0,1]} x^m (1-x)^{n-m} \, dF(x), \qquad 0 \le m \le n.$$

This probabilistic interpretation of Hausdorff's theorem can be further developed as follows.

**Theorem 1.3** (B. de Finetti's theorem [9], §VII.4)**.** *For $c$ completely monotone with $c(0) = 1$, the triangular array (1.4) defines the probability distribution of the random number $S_n$ of successes in the first $n$ trials of an infinite exchangeable sequence of successes and failures, according to (1.6). For such an exchangeable sequence, $S_n/n$ converges almost surely to a random variable $X$ with distribution function $F$, and conditionally given $X$, the trials are independent with success probability $X$.*

In a recent study of moment problems [7], Diaconis and Freedman considered the family of probability distributions $F$ on $[0,1]$ with a nondecreasing density $f$ on the semi-open interval $[0,1[$, so that

$$F(x) = \int_{[0,x]} f(u) du, \qquad 0 \le x < 1.$$

Note that $F(1-) < 1$ is allowed, in which case the distribution has an atom of magnitude $1 - F(1-)$ at 1. Equivalently, $F(x)$ is a convex function of $x \in [0,1]$, with $F(0) = 0$. Diaconis and Freedman showed the following.

**Theorem 1.4** ([7], Theorem 10)**.** *A sequence $c$ with $c(0) = 1$ is the moment sequence of a convex probability distribution function $F$ on $[0,1]$ with $F(0) = 0$ if and only if for each fixed $n = 1, 2, \ldots$ the sequence $c(n,m)$ is non-negative and nondecreasing in $m$:*

$$(1.9) \qquad c(n,0) \ge 0 \quad \text{and} \quad c(n,m+1) - c(n,m) \ge 0, \qquad 0 \le m < n.$$

In view of (1.6) and (1.7), condition (1.9) is a natural discrete analogue of the nondecreasing density condition for $F$.

In the course of our work on partition structures derived from regenerative random sets [12, 13, 14], we faced a similar problem of characterizing the Laplace exponent of a subordinator in terms of its values at positive integers. We noticed that this problem was equivalent to the moment problem for convex distributions on $[0,1]$, and that both problems can be reduced to a known integral representation of sequences $a$ subject to the following condition, which was studied by Choquet. See [1, Definition 6.1, p. 130] for history and terminology.

**Definition 1.5.** A sequence $a$ is *completely alternating* if the sequence $-\nabla a$ is completely monotone. That is to say:

$$(1.10) \qquad \nabla^j a(n) \leq 0, \quad j = 1, 2, \ldots, \; n = 0, 1, \ldots,$$

The fundamental representation of completely alternating sequences is provided by the following theorem.

**Theorem 1.6** (Contained in [1], Proposition 6.12). *A sequence $a$ is completely alternating if and only if there exists a bounded measure $\nu$ on $[0,1[$ and constants $A \geq 0$ and $B$ such that*

$$(1.11) \qquad a(n) = An + B + \int_{[0,1[} \frac{1-x^n}{1-x} \nu(dx).$$

*The data $(A, B, \nu)$ are uniquely determined.*

In Section 2 we show how Theorem 1.6 follows easily from Hausdorff's Theorem 1.2, and then deduce the following variant of Theorem 1.4.

**Theorem 1.7.** *A sequence $c$ is the sequence of moments of a convex distribution function on $[0,1]$ with $F(0) = 0$ if and only if $c(0) = 1$ and the sequence $a$ defined by*

$$(1.12) \qquad a(0) = 0, \quad a(n) = n\,c(n-1), \quad n = 1, 2, \ldots$$

*is completely alternating.*

Comparing Theorems 1.4 and 1.7, we deduce that condition (1.9) on a sequence $c$ with $c(0) = 1$ is equivalent to the condition that $a$ derived from $c$ via (1.12) is completely alternating. We check this directly by algebra in Section 3, thereby providing a new proof of Theorem 1.4. We explain in Section 4 how we were first led to Theorem 1.7 by consideration of the Lévy-Khintchine formula for the Laplace transform of a subordinator. Finally, Section 6 relates distributions with higher convexity properties to alternating sequences of higher order.

## 2. Proofs of Theorems 1.6 and 1.7

*Proof of Theorem 1.6.* Let $a$ be a completely alternating sequence with $a(0) = 0$. Then $-\nabla a$ is a completely monotone sequence which can be represented by Hausdorff's theorem as

$$(2.1) \qquad -\nabla a(n) = \int_{[0,1]} \xi^n \nu(d\xi), \quad n = 0, 1, \ldots$$

for a unique bounded measure $\nu$ on $[0,1]$. Hence, by summing a geometric series

$$(2.2) \qquad a(n) = n\nu\{1\} + \int_{[0,1[} \frac{1-\xi^n}{1-\xi} \nu(d\xi), \quad n = 0, 1, \ldots.$$

This is the special case of (1.11) with $B = 0$. This case, applied after first subtracting $a(0)$ from all terms of a completely alternating sequence $a$, gives the general form (1.11) with $B = a(0)$. Conversely, if $a$ is defined by (1.11) or its special case (2.2), then (2.1) is obtained by subtraction, hence $-\nabla a$ is the completely monotone moment sequence of $\nu$. □

We note in passing that the linear term $n\nu\{1\}$ in (2.5) could be absorbed into the integral by extending the integral from $[0,1[$ to $[0,1]$, with evaluation of the integrand by continuity at $\xi = 1$. But in this and the similar expressions (1.11) and (2.2) we prefer to display this term separately to avoid any possible misunderstanding.

Our proof of Theorem 1.7 is based on the integral representation of probability distributions with convex distribution functions, used also by Diaconis and Freedman. Let

$$F_\xi(x) := \frac{x-\xi}{1-\xi} \mathbf{1}(x \geq \xi), \quad \xi \in [0,1[,$$

which is the distribution function of the uniform distribution on $[\xi, 1]$, with density

$$f_\xi(x) = \frac{1}{1-\xi} \mathbf{1}(x \geq \xi), \quad x \in [0,1].$$

Let $F_1(x) = \mathbf{1}(x=1)$, corresponding to a unit mass at 1.

**Lemma 2.1** ([7], Lemma 2). *The formula*

$$(2.3) \qquad F(x) = \int_{[0,1]} F_\xi(x) \nu(\mathrm{d}\xi),$$

*sets up a bijection between convex probability distribution functions $F$ on $[0,1]$ with $F(0) = 0$ and probability measures $\nu$ on $[0,1]$. This relation between $F$ and $\nu$ implies*

$$(2.4) \qquad \int_{[0,1]} g(x) \, dF(x) = \int_{[0,1]} \left( \int_{[0,1]} g(x) \, dF_\xi(x) \right) \nu(d\xi)$$

*for every non-negative measurable $g$.*

*Proof.* A convex distribution function $F$ with $F(0) = 0$ has on $[0,1[$ a nondecreasing density $f$, a version of which is the right derivative of $F$. Hence it is clear that a unique probability measure $\nu$ is defined by

$$f(x) = \int_{[0,x]} \frac{\nu(d\xi)}{1-\xi}, \quad x \in [0,1[,$$

and $\nu\{1\} = F(1) - F(1-)$, from which (2.4) follows by Fubini's theorem. □

**Lemma 2.2.** *Let $F$ be the mixture of $F_\xi$'s with respect to $\nu(d\xi)$ for some probability measure $\nu$ on $[0,1]$, as in Lemma 2.1. Then the moments of $F$ are determined by*

$$(2.5) \qquad n \int_{[0,1[} x^{n-1} \, dF(x) = n\nu\{1\} + \int_{[0,1[} \frac{1-\xi^n}{1-\xi} \nu(d\xi), \quad n = 1, 2, \ldots.$$

*Let $a(0) = 0$ and let $a(n)$ be the common value of both sides of (2.5) for $n = 1, 2, \ldots$. Then $a$ is completely alternating. Indeed $-\nabla a$ is the completely monotone moment sequence of $\nu$, as in (2.1).*

*Proof.* By elementary integration, the moments of $F_\xi$ are given by

$$(2.6) \qquad n \int_{[0,1]} x^{n-1} \, dF_\xi(x) = \frac{1 - \xi^n}{1 - \xi}, \quad n = 1, 2, \ldots.$$

So (2.5) is the instance of (2.4) with $g(x) = nx^{n-1}$. The rest is read from the previous discussion of (2.2) and (2.1), as formalized in Theorem 1.6. □

*Proof of Theorem 1.7.* If $c$ is the moment sequence of $F$, and $a$ is derived from $c$ via (1.12), then $a$ is completely alternating by the previous lemma. Conversely, given a sequence $c$ with $c(0) = 1$, let $a$ be derived from $c$ via (1.12). If $a$ is completely alternating, then Theorem 1.6 represents $a(n)$ by the right side of (2.5) for $n = 1, 2, \ldots$, for some unique probability measure $\nu$ on $[0, 1]$. Then, by (2.5), the sequence of moments of $F$ is the same as if $F$ were the convex distribution function uniquely associated with $\nu$ via (2.4). Finally, $F$ equals this convex distribution function, by uniqueness of the solution of the Hausdorff moment problem. □

**Remark.** A subtle point of the above argument is that possibly $\int_{[0,1[} x^{-1} \, dF(x) = \infty$, in which case the left side of (2.5) has no meaning for $n = 0$. We insist in any case that $a(0) = 0$ in Theorem 1.7 and Lemma 2.2 *by definition*. Later discussion in Section 4 makes it clear that 0 is the limiting value of $\lambda \int_{[0,1[} x^{\lambda-1} \, dF(x)$ as $\lambda \downarrow 0$. But this fact is not relevant to the present argument.

## 3. Some algebra

Here we check directly that condition (1.9) on a sequence $c$ with $c(0) = 1$ is equivalent to the condition that $a$ derived from $c$ via (1.12) is completely alternating. According to the Leibnitz rule from the calculus of finite differences, the product $xy$ of two sequences $x$ and $y$ has successive differences

$$(3.1) \qquad \nabla^j(xy)(n) = \sum_{i=0}^{j} \binom{j}{i} \nabla^{j-i} x(n+i) \nabla^i y(n).$$

Applied to $a(n) := nc(n-1)$, with $c(-1) := 0$, this gives

$$(3.2) \qquad \nabla^j a(n) = n \nabla^j c(n-1) - j \nabla^{j-1} c(n).$$

A simple computation using (1.4) and (3.2) now shows that

$$(3.3) \qquad c(n, m+1) - c(n, m) = \frac{-1}{m+1} \binom{n}{m} \nabla^{n-m} a(m+1), \quad 0 \le m < n.$$

The equivalence of (1.9) and the completely alternating condition on $a$ is now evident from (3.2) for $n = 0$ and (3.3).

## 4. Subordinators

We explain in this section how we were led to formulate Theorem 1.7 by the appearance of completely alternating sequences and monotone densities in another context: the Lévy-Khintchine formula for the Laplace transform of a subordinator.

Let $H$ denote a probability distribution on $[0, \infty]$. Recall that $H$ is called *infinitely divisible* if for every $n = 2, 3, \ldots$ there exists a sequence of $n$ independent and identically distributed random variables, with values in $[0, \infty]$, whose sum is distributed according to $H$. It is well known that $H$ is infinitely divisible with $H[0, \infty[ > 0$ if and only if $H$ is the distribution of $Y(1)$ for some *subordinator* $(Y(t),\ t \geq 0)$, that is a process with stationary independent nonnegative increments with values in $[0, \infty]$, and $Y(0) = 0$. Note that this process may jump to the terminal value $\infty$ in finite time, but is not identically $\infty$.

**Theorem 4.1** (Lévy-Khintchine formula [2])**.** *The Laplace transform of a subordinator* $(Y(t),\ t \geq 0)$ *is given by the formula*

$$(4.1) \qquad \mathbb{E}\left[e^{-\lambda Y(t)}\right] = e^{-t\Phi(\lambda)}$$

*with Laplace exponent*

$$(4.2) \qquad \Phi(\lambda) = \lambda \mathtt{d} + \int_{]0,\infty]} (1 - e^{-\lambda y}) \Lambda(dy)$$

*for some* $\mathtt{d} \geq 0$ *and some measure* $\Lambda$ *on* $]0, \infty]$, *with* $(\mathtt{d}, \Lambda)$ *determined uniquely by the distribution of* $Y(t)$ *for any fixed* $t > 0$.

**Remark.** The constant $\mathtt{d}$ accounts for the continuous drift component of the subordinator, while $\Lambda$, called the *Lévy measure*, governs the Poisson rates of jumps of various sizes. In particular, the mass $\Lambda\{\infty\}$ is the rate at which the process jumps to $\infty$.

*Proof.* We sketch only the following derivation of (4.1) and (4.2), adapted from [2, pp. 5-7], because it highlights the connection with monotone density problems. Let $(Y(t),\ t \geq 0)$ be a subordinator. By independence and stationarity of increments

$$\mathbb{E}\left[e^{-\lambda Y(s+t)}\right] = \mathbb{E}\left[e^{-\lambda Y(s)}\right] \mathbb{E}\left[e^{-\lambda Y(t)}\right],$$

which implies (4.1) with $0 < \Phi(\lambda) < \infty$ by virtue of the well known functional equation for the exponential function. Using again the properties of increments,

$$\Phi(\lambda) = \lim_{m \to \infty} m(1 - e^{-\Phi(\lambda)/m}) \;=\; \lim_{m \to \infty} m \mathbb{E}[1 - e^{-\lambda Y(1/m)}]$$

$$= \lambda \lim_{m \to \infty} \int_{]0,\infty[} e^{-\lambda y}\, m\, \mathbb{P}(Y(1/m) \geq y)\, dy\,,$$

where $m\, \mathbb{P}(Y(1/m) \geq y)\, dy$ is a measure with nonincreasing density for each $m = 1, 2, \ldots$. By continuity properties of the Laplace transform these measures must converge, as $m \to \infty$, to a measure $\overline{\Lambda}(y)dy$ with nonincreasing density and a possible atom $\mathtt{d}$ at $0$. Thus

$$(4.3) \qquad \Phi(\lambda) = \lambda \mathtt{d} + \lambda \int_{]0,\infty[} e^{-\lambda y} \overline{\Lambda}(y) dy$$

from which (4.2) follows by integration by parts, with measure $\Lambda$ defined through $\Lambda\,]y, \infty] = \overline{\Lambda}(y)$, $y \geq 0$, for a right-continuous choice of the monotone density $\overline{\Lambda}(y)$. □

To relate the Lévy-Khintchine formula to the discussion of the previous sections, make the change of variables $x = e^{-y}$ in (4.3) to rewrite this equation as

$$(4.4) \qquad \Phi(\lambda) = \lambda \mathtt{d} + \lambda \int_{]0,1[} x^{\lambda-1} \overline{\Lambda}(-\log x)\, dx\,,$$

or again

$$(4.5) \qquad \Phi(\lambda) = \lambda \int_{]0,1]} x^{\lambda-1} dF(x),$$

where $F(x)$ is the convex distribution function with $F(0) = 0$, with increasing derivative $\overline{\Lambda}(-\log x)$ for $x \in [0, 1[$ and with an atom of magnitude d at 1. Formula (4.5) shows that the problem of characterizing moment sequences of convex probability distribution functions on $[0, 1]$ with $F(0) = 0$ is identical to the problem of characterizing the sequence of evaluations $\Phi(1), \Phi(2), \ldots$ for the Laplace exponent $\Phi$ of a subordinator subject to the normalization condition $\Phi(1) = 1$. Since (4.2) can be rewritten as

$$(4.6) \qquad \Phi(\lambda) = \lambda \mathtt{d} + \int_{[0,1]} \frac{1-x^\lambda}{1-x} \nu(dx),$$

where $\nu(\mathrm{d}x)$ is the image of $(1 - e^{-z})\Lambda(\mathrm{d}z)$ via $x = e^{-z}$, we deduce the following corollary of Theorem 1.6.

**Corollary 4.2.** *There exists a subordinator with Laplace exponent $\Phi$ subject to a prescribed sequence of values*

$$\Phi(0) = 0, \ \Phi(1), \ \Phi(2), \ldots$$

*if and only if this sequence is completely alternating. The value of $\Phi(\lambda)$ is then determined for all $\lambda \geq 0$ by the Lévy-Khintchine formula (4.2), for some uniquely determined d and $\Lambda$.*

Theorem 1.7 can be read from Corollary 4.2 and the alternate representation (4.5) of the Laplace exponent. This is how we first recognized that the completely alternating condition of Theorem 1.7 was equivalent to the condition (1.9) in Theorem 1.4.

**Further remarks.** According to a classical theorem of Münz [18], a sufficient condition for density in $C[0, 1]$ of the space of functions spanned on $x \mapsto x^{\lambda_n}$, $n = 1, 2, \ldots$, is that $\lambda_n \to \infty$ and $\sum_n 1/\lambda_n = \infty$. So the value of $\Phi(\lambda)$ at every real $\lambda \geq 0$ is determined by the values $\Phi(\lambda_n)$ at points of such a sequence $(\lambda_n)$. Newton's interpolation series [10] allows to explicitly recover $\Phi(\lambda)$ from the discrete evaluations $\Phi(\lambda_n)$, $n = 1, 2, \ldots$. In the same vein as Corollary 4.2, Hausdorff's theorem implies that a sequence $\Psi(0) = 1, \Psi(1), \Psi(2), \ldots$ is the sequence of evaluations of

$$\Psi(\lambda) = \mathbb{E}(e^{-\lambda X})$$

for some distribution of $X$ on $[0, \infty]$ if and only if this sequence is completely monotone. In that case the distribution of $X$, and hence $\Psi(\lambda)$ for all $\lambda \geq 0$, is uniquely determined. As observed by Feller [8], consideration of the sequences $\Psi(0) = 1, \Psi(\epsilon), \Psi(2\epsilon), \ldots$ for $\epsilon = 2^{-k}$ and $k = 0, 1, 2, \ldots$, leads to *Bernstein's theorem* that a function $\Psi$ is the Laplace transform of some probability distribution on $[0, \infty[$ if and only if $\Psi$ is infinitely differentiable with $\Psi(0) = 1$ and $(-1)^n \Psi^{(n)} \geq 0$ for every $n$, where $\Psi^{(n)}$ is the $n$th derivative of $\Psi$.

Similarly, a function $\Phi$ is the Laplace exponent of some subordinator if and only if $\Phi$ is infinitely differentiable with $\Phi(0) = 0$ and $(-1)^n \Phi^{(n)} \leq 0$ for every $n$. See [9] for further discussion, and [1] for the theory of completely monotone and completely alternating functions defined on a semigroup instead of the positive integers or the positive halfline.

## 5. Regenerative composition structures

We sketch in this section a probabilistic interpretation of completely alternating sequences, based on our recent work [12, 13, 14] on regenerative compositions and their associated partition structures. We note in passing that developments and applications of de Finetti's theorem have played an important role in a number of studies of measure-valued processes and associated particle systems [4, 5, 6]. The general idea is that some random allocation or splitting of a mass continuum can be described in terms of a simpler combinatorial model obtained by independent random sampling of $n$ points in the continuum. Probabilities in the combinatorial model are typically represented as moments of random variables of interest in the continuum model. The combinatorial models are consistent in a natural sense as the sample size $n$ varies, and the continuum model is recovered as a law of large numbers limit.

The simplest illustration of this idea is provided by de Finetti's theorem for sequences of zeros and ones (Theorem 1.3). Suppose the unit interval is split by a point $X$ into two interval components $[0, X]$ and $]X, 1]$. Let $U_1, U_2, \ldots$ be a sequence of independent uniform $[0, 1]$ variables independent of $X$, and let

$$S_n = \sum_{i=1}^{n} 1(U_i \leq X).$$

Associating each $U_i$ with a 'ball', and $[0, X]$ and $]X, 1]$ with two 'boxes', $S_n$ describes an allocation of $n$ balls in two boxes, with $S_n$ balls in the left-hand box, and $n - S_n$ balls in the right-hand box. As $n$ varies, these allocations are *sampling consistent*, meaning that if a ball is picked uniformly at random and deleted from the $n$th random allocation, the result is distributed like the $(n-1)$th random allocation. To paraphrase de Finetti's theorem: every sampling consistent sequence of distributions for the allocation of balls $S_n$ can be realized via this scheme directed by some random variable $X \in [0, 1]$, and $X$ is recovered as the limit of $S_n/n$ as $n \to \infty$.

A straightforward generalization of this model provides an interpretation of de Finetti's theorem for exchangeable sequences with a finite or countably infinite number of possible values. Let

$$0 = X_0 \leq X_1 \leq X_2 \leq \cdots \leq 1$$

with $\lim_n X_n = 1$ be the cumulative random sums

$$X_k = \sum_{j=1}^{k} P_k$$

associated with some *random discrete distribution* $(P_1, P_2, \ldots)$. Let $U_1, U_2, \ldots$ be a sequence of independent uniform $[0, 1]$ variables independent of the $X_k$'s, and set

$$N_{n,j} := \sum_{i=1}^{n} 1(X_{j-1} < U_i \leq X_j), \quad j = 1, 2, \ldots.$$

Then the random sequence of counts

$$(N_{n,1}, N_{n,2}, \ldots) \quad \text{with} \quad \sum_{j=1}^{\infty} N_{n,j} = n$$

can be interpreted as a random allocation of $n$ balls into a sequence of boxes labeled by $j = 1, 2, \ldots$. These allocations are sampling consistent in an obvious sense, and according to de Finetti's theorem every sampling consistent sequence of allocations can be represented this way, with $P_j$ the limiting frequency of balls in box $j$.

In some applications of random discrete distributions, the labeling of atoms of the distribution is of little or no importance. If all that matters is the relative sizes of these atoms, and perhaps some ordering of these atoms (not necessarily a simple indexing by positive integers), the natural combinatorial object is an allocation of $n$ balls into some finite number of non-empty boxes, which might be either ordered or unordered. For $n$ a positive integer, a *composition of $n$* is a sequence of positive integers with sum $n$. A *partition of $n$* is a non-increasing composition of $n$. The terms of a composition may be called its *parts*. A *composition structure* is a sequence $(\mathcal{C}_n, n = 1, 2, \ldots)$ of random compositions of integers which is *sampling consistent*, that is:

> If $n$ identical balls are distributed into an ordered series of boxes according to $(\mathcal{C}_n)$, then a distributional copy of $\mathcal{C}_{n-1}$ is obtained by discarding one of the balls picked uniformly at random, independently of $(\mathcal{C}_n)$, and then deleting an empty box in case one is created.

This is a variation of Kingman's notion of a *partition structure* [17], which is a sequence of random partitions of integers $(\mathcal{P}_n)$ subject to the same consistency condition, except that after discarding a ball and deleting an empty box if necessary, the boxes are permuted to obtain a partition of $n - 1$. Kingman [17] established a one-to-one correspondence between partition structures $(\mathcal{P}_n)$ and distributions for a sequence of nonnegative random variables $V_1, V_2, \ldots$ with $V_1 \geq V_2 \geq \ldots$ and $\sum_i V_i \leq 1$. In Kingman's *paintbox representation*, the random partition $\mathcal{P}_n$ of $n$ is constructed as follows from $(V_k)$ and a sequence of independent random variables $U_i$ with uniform distribution on $[0, 1]$, where $(U_i)$ and $(V_k)$ are independent: $\mathcal{P}_n$ is defined to be the sequence of ranked sizes of blocks of the partition of $[n]$ generated by a random equivalence relation $\sim$ on positive integers, with $i \sim j$ if and only if either $i = j$ or both $U_i$ and $U_j$ fall in $I_k$ for some $k$, where the $I_k$ are some disjoint random sub-intervals of $[0, 1]$ of lengths $V_k$. See also [20] and papers cited there for further background.

Gnedin [11] gave a similar representation of composition structures, using a random closed $\mathcal{R} \subset [0, 1]$ to separate points of a uniform sample into clusters. Given $\mathcal{R}$, define an interval partition of $[0, 1]$ comprised of *gaps*, that is open interval components of $[0, 1] \setminus \mathcal{R}$, and of individual points of $\mathcal{R}$. A random ordered partition of the set $[n] := \{1, \ldots, n\}$ is constructed from $\mathcal{R}$ and independent uniform sample points $U_1, \ldots, U_n$ by grouping the indices of sample points which fall in the same gap, and letting the points which hit $\mathcal{R}$ to be singletons. A random composition $\mathcal{C}_n$ of $n$ is then constructed as the sequence of block sizes in this partition of $[n]$, ordering the blocks from left to right, according to the location of the corresponding sample points in $[0, 1]$. Gnedin showed that every composition structure $(\mathcal{C}_n)$ can be so represented. As in Kingman's representation of partition structures, $\mathcal{R}$ can be interpreted as an asymptotic shape of $\mathcal{C}_n$, provided $\mathcal{C}_n$ is properly encoded as an element of the metric space of closed subsets of $[0, 1]$ with the Hausdorff distance function.

An interesting class of composition structures $(\mathcal{C}_n)$ is obtained by supposing that $\mathcal{R}$ is the closure of $\{1 - \exp(-Y(t)), t \geq 0\}$ for some subordinator $(Y(t), t \geq 0)$. As shown in [12], these composition structures are characterized by the following *regenerative property*:

For all $n > m \geq 1$, given that the first part of $\mathcal{C}_n$ is $m$, the remaining composition of $n - m$ is distributed like $\mathcal{C}_{n-m}$.

Let $q(n, m)$ denote the probability that the first part of $(\mathcal{C}_n)$ is of size $m$. Then for each composition $\lambda = (\lambda_1, \ldots, \lambda_\ell)$ of $n$

$$(5.1) \qquad \mathbb{P}(\mathcal{C}_n = \lambda) = \prod_{j=1}^{\ell} q(\Lambda_j, \lambda_j),$$

where $\Lambda_j = \lambda_j + \ldots + \lambda_\ell$. We showed in [12] that this formula together with a recursion implied by sampling consistency forces

$$(5.2) \qquad q(n, m) = \frac{-\binom{n}{m} \nabla^m \Phi(n - m)}{\Phi(n)}$$

for some unique sequence $\Phi(n)$, $n = 0, 1, 2, \ldots$ with $\Phi(0) = 0$ and $\Phi(1) = 1$. Non-negativity of the matrix $q$ shows that $\Phi$ is completely alternating. So Corollary 4.2 now provides the integral representation

$$(5.3) \qquad \Phi(n) = n\mathtt{d} + \int_{]0,1]} (1 - (1 - x)^n) \tilde{\nu}(dx).$$

for $\tilde{\nu}(dx)$ the image of the Lévy measure $\Lambda(\mathrm{d}y)$ of the subordinator via the transformation from $y$ to $1 - \exp(-y)$, and hence

$$(5.4) \qquad q(n, m) = n\mathtt{d}\, 1(m = 1) + \binom{n}{m} \int_{]0,1]} x^m (1 - x)^{n-m} \tilde{\nu}(dx).$$

## 6. Higher convexity

A distribution function $F$ on $[0, 1]$ is said to be $(k + 1)$-*convex* if the derivative $F^{(k-1)}$ exists and is a convex function on $[0, 1[$ (see [3, 16, 19] and references therein for this and other generalized concepts of convexity and their applications in statistics). The moment problem for this class of distributions can be analyzed in the same way as for the convex distributions we considered above (which were 2-convex, with $k = 1$). For a $(k + 1)$-convex $F$, Taylor's formula with remainder becomes

$$(6.1) \qquad F(x) = P(x) + \int_{[0,x]} \frac{(x - \xi)^k}{(1 - \xi)^k} \nu(d\xi), \quad x \in [0, 1[,$$

where $\nu$ is a bounded measure on $[0, 1]$ and $P$ is a polynomial of degree at most $k - 1$. We suppose further that $P = 0$. In this case the analogue of (2.5) becomes

$$(n+1)^{k\uparrow} c(n) = A\,(n+1)^{k\uparrow} + k! \int_{[0,1[} \left[ x^{n+k} - \sum_{j=0}^{k-1} \binom{n+k}{j}(x-1)^j \right] \frac{\nu(dx)}{(x-1)^k},$$

where $n^{k\uparrow} = n(n+1) \cdots (n+k-1)$ and $A \geq 0$. A principal role is delegated to the $k$-*associated* sequence

$$a(0) = \ldots = a_{k-1} = 0, \quad a_k(n) = n^{k\downarrow}\, c(n - k) \quad \text{for } n \geq k,$$

which is $k$-*alternating* sequence meaning that $(-\nabla)^k a_k$ is completely monotone. A generalization of Proposition 1.7 emerges.

**Proposition 6.1.** *For each $k = 1, 2, \ldots$, a distribution function on $[0,1]$ is $(k+1)$-convex with $F^{(j)}(0+) = 0$ for $j = 0, \ldots, k-1$ if and only if $c(0) = 1$ and the $k$-associated sequence $a_k$ is $k$-alternating.*

The side condition in the proposition kills the polynomial part $P$ in (6.1), a condition which was guaranteed by $F(0) = 0$ in the case $k = 1$. If the coefficients of $P$ are nonnegative, it is still possible to give a characterization in terms of a $k$-associated sequence alone, whose initial $k$ terms are no longer zeroes, rather need to be defined in a nontrivial way by Newton's interpolation.